\documentclass[a4paper, 12pt]{article}
\usepackage[utf8]{inputenc}
\usepackage{amssymb,amsmath,amsfonts}
\usepackage{latexsym}
\usepackage{mathrsfs}
\usepackage{amsthm}
\usepackage[left=2cm, right=2cm, top=3cm, bottom=3cm]{geometry}
\usepackage{graphicx}
\usepackage{enumerate}
\usepackage{tikz}
\usepackage{tabularx}

\usetikzlibrary{decorations,arrows}
\usetikzlibrary{decorations.pathmorphing}
\usepgflibrary{decorations.pathreplacing}

\newcommand{\To}{\rightarrow}

\newcommand{\h}{\mathcal{H}(K_{s,t,\alpha,\beta,\bsgamma})}

\newcommand{\bszero}{\boldsymbol{0}}
\newcommand{\bsgamma}{\boldsymbol{\gamma}}

\newcommand{\boldg}{\boldsymbol{g}}
\newcommand{\bsh}{\boldsymbol{h}}

\newcommand{\bsk}{\boldsymbol{k}}
\newcommand{\bsl}{\boldsymbol{l}}
\newcommand{\bsm}{\boldsymbol{m}}

\newcommand{\boldp}{\boldsymbol{p}}
\newcommand{\boldq}{\boldsymbol{q}}

\newcommand{\bsx}{\boldsymbol{x}}
\newcommand{\bsy}{\boldsymbol{y}}
\newcommand{\boldz}{\boldsymbol{z}}
\newcommand{\uprighte}{\mathrm{e}}
\newcommand{\NN}{\mathbb{N}}
\newcommand{\ZZ}{\mathbb{Z}}
\newcommand{\FF}{\mathbb{F}}

\newcommand{\LL}{\mathbb{L}}
\newcommand{\wal}{\mathrm{wal}}
\newcommand{\walarg}[2]{\mathrm{wal}_{#1{}}(#2{})}
\newcommand{\trig}[1]{\uprighte_{#1{}}}
\newcommand{\trigarg}[2]{\uprighte_{#1{}}(#2{})}
\newcommand{\ii}{\mathtt{i}}

\newcommand{\linalg}[1]{A_{N,s,t,M}{(f)(#1{})}}

\newcommand{\EMB}{\mathrm{EMB}}

\newcommand{\cH}{{\cal H}}

\newtheorem{thm}{Theorem}
\newtheorem{cor}{Corollary}
\newtheorem{prop}{Proposition}
\newtheorem{lem}{Lemma}

\newtheorem{rem}{Remark}

\allowdisplaybreaks

\begin{document}
\title{Tractability of $\mathbb{L}_2$-approximation in hybrid function spaces}
\author{Peter Kritzer\thanks{P. Kritzer is supported by the Austrian
Science Fund (FWF): Project F5506-N26.}\;,
Helene Laimer\thanks{H. Laimer is supported by the Austrian 
Science Fund (FWF): Projects F5506-N26 and F5508-N26.}\;,
Friedrich Pillichshammer\thanks{F. Pillichshammer is supported by the Austrian 
Science Fund (FWF): Project F5509-N26.\newline All projects are parts of the 
Special Research Program "Quasi-Monte Carlo Methods: Theory and Applications".}
}
\date{}

\maketitle

\begin{abstract}
We consider multivariate $\mathbb{L}_2$-approximation in reproducing kernel Hilbert
spaces which are tensor products of weighted Walsh spaces and weighted Korobov
spaces. We study the minimal worst-case error $e^{\mathbb{L}_2-\mathrm{app},\Lambda}(N,d)$ of all algorithms that use $N$ information evaluations from the class $\Lambda$ in the $d$-dimensional case.   
The two classes $\Lambda$ considered in this paper are the class $\Lambda^{{\rm all}}$ consisting of all linear functionals and the class $\Lambda^{{\rm std}}$ consisting only of function evaluations.

The focus lies on the dependence of $e^{\mathbb{L}_2-\mathrm{app},\Lambda}(N,d)$ on the dimension $d$. The main results are conditions for weak, polynomial, and strong polynomial tractability.
\end{abstract}



\section{Introduction}

We consider $\LL_2$-approximation of functions in certain reproducing kernel Hilbert spaces $\cH(K)$,
which are embedded into $\LL_2 ([0,1]^d)$, where $K$ denotes the reproducing kernel. To be more precise, we approximate the embedding operator  $$\EMB_{d}:\cH(K)\rightarrow \LL_2([0,1]^d), \ \ \EMB_{d}(f)=f,$$ and measure the approximation error in the $\LL_2$-norm. Since $\cH(K)$ is a reproducing kernel Hilbert space it is known (cf. \cite{NW08,TWW88}) that there is no loss of generality when we restrict ourselves
to linear approximation algorithms of the form $A_{N,d}(f) = \sum_{k=1}^{N}{a_k L_k (f)}$
with coefficients $a_k \in \LL_2([0,1]^d)$ and continuous linear functionals $L_k$ on $\cH (K)$ from a permissible class of information $\Lambda$. Here $N$ is the number of information evaluations. 

We study the problem in the so-called worst-case setting, i.e., we measure the approximation error of an algorithm $A_{N,d}$ by means of the worst-case error, 
\begin{align*}
e^{\LL_2-\mathrm{app}}(A_{N,d}) = 
\sup_{\substack{f \in \cH (K) \\ \|f\|_{\cH(K)}\leq 1}}{\|\EMB_{d}{(f)} - A_{N,d}(f)\|_{\LL_2([0,1]^d)}}.
\end{align*}
The $N^{{\rm th}}$ minimal worst-case error is given by
$$
e^{\LL_2-\mathrm{app},\Lambda}(N,d)=\inf_{A_{N,d}}e^{\mathrm{app}}(A_{N,d}),
$$
where the infimum is extended over all linear algorithms $A_{N,d}$ using $N$ information evaluations 
from the class $\Lambda$. We are particularly interested in the dependence of the $N^{{\rm th}}$ minimal worst-case 
error on the dimension $d$. To study this dependence systematically we consider the information complexity $N^{\LL_2-\mathrm{app},\Lambda}(\varepsilon,d)$, which is the minimal number $N$ for which 
there exists an algorithm using $N$ information evaluations from the class $\Lambda\in\{\Lambda^{{\rm all}},\Lambda^{{\rm std}}\}$ with an error of at most $\varepsilon$.

We would like to avoid cases where the information complexity $N^{\LL_2-\mathrm{app},\Lambda}(\varepsilon,d)$ grows exponentially or even faster with the dimension $d$ or with $\varepsilon^{-1}$. To 
quantify the behavior of the information complexity we use different notions of tractability, namely strong polynomial tractability, polynomial tractability, and weak tractability (we refer to Section 
\ref{secproblem} for the precise definitions).

The current state of the art of tractability theory is summarized in the three volumes of the book of Novak and Wo\'{z}niakowski \cite{NW08,NW10,NW12} 
which we refer to for extensive information on this subject and further references.

In previous papers, several authors have studied similar approximation problems 
in various reproducing kernel Hilbert spaces, see, e.g., \cite{BDK09, DKK08, DKPW14, KSW06, NSW,ZKH09}. These investigations have in common that the  
reproducing kernel Hilbert spaces considered 
are tensor products of one-dimensional spaces whose kernels are all of the same type (but maybe equipped with different weights). 
In the current paper we consider the case where the reproducing kernel is a product of kernels of different type. 
We call such spaces {\it hybrid spaces}. Some results on tractability in general hybrid spaces can be found in the literature. For example, 
in \cite{NW10} multivariate integration is studied for arbitrary reproducing kernels $K_d$ without relation to $K_{d+1}$.  
Here we consider as a special instance the tensor product of Walsh and Korobov spaces. The problem of numerical integration in such spaces was recently 
considered in \cite{KP15}. The study of a hybrid 
of Korobov and Walsh spaces could be important in view of 
functions which are periodic with respect to some of the components and, for example, piece-wise constant with respect to the remaining components. 
Moreover, it has been pointed out by several scientists (see, e.g., \cite{K13,L14}) that hybrid problems may be relevant for 
certain applications. Indeed, communication with the authors of \cite{K13} and \cite{L14} have motivated our idea 
for considering function spaces where we may have very different properties of the elements with respect to different components, as for 
example regarding smoothness. 

From the analytical point of view, it is very challenging to deal with hybrid spaces. 
The reason for this is the rather complex interplay between the different analytic and algebraic structures of the kernel functions. 
In the present study we are concerned with Fourier analysis carried out simultaneously with respect to the Walsh and the trigonometric function systems. 
The problem is also closely related to the study of hybrid point sets which received much attention in recent years 
(see, for example, \cite{H11,HKLP09}). Hence we also have considerable theoretical interest in studying this problem.

The paper is organized as follows. In Section \ref{secspace} we introduce the Hilbert space under consideration. In Section \ref{secproblem} 
we specify the problem setting and state our main result. The proofs are presented in Section \ref{sec:proof}.

\section{The hybrid function space}\label{secspace}

We study a specific reproducing kernel Hilbert
space, namely the tensor product of a Korobov space and a Walsh
space, that was introduced in \cite{KP15}. See \cite{A50} for general information about reproducing kernel Hilbert spaces.

Fix a prime number $b$ and let $\ii=\sqrt{-1}$. For $k \in \mathbb{N}_0$ with $b$-adic expansion
$k = \kappa_a b^a +\cdots +\kappa_1 b+\kappa_0$ with $\kappa_j \in \{0, \dots , b-1 \}$ we define the $k^{{\rm th}}$ Walsh function $\wal_k : [0,1) \rightarrow \mathbb{C}$ by
\begin{align*}
\wal_k(x) = \exp\left( 2 \pi \ii \frac{\xi_1 \kappa_0 + \dots + \xi_{a+1} \kappa_a}{b} \right),
\end{align*}
for $x \in [0,1)$ with $b$-adic expansion $x=\frac{\xi_1}{b}+\frac{\xi_2}{b^2}+\cdots$ (unique in the sense that infinitely many of the $\xi_i$ are different from $b-1$). Note that $a = \left\lfloor 
\log_b{k} \right\rfloor$. 

For $\bsk = (k_1, \dots ,k_s) \in \mathbb{N}_0^s$ and $\bsx =
(x_1, \dots , x_s) \in [0,1)^s$ the $\bsk^{{\rm th}}$ $s$-variate Walsh function 
$\wal_{\bsk} : [0,1)^s \rightarrow \mathbb{C}$
is given by
$\wal_{\bsk}(\bsx) = \prod_{j=1}^s \wal_{k_j}(x_j)$.

Further, for $\bsl \in \mathbb{Z}^t$ we define the $t$-variate $\bsl^{{\rm th}}$ trigonometric function 
$\trig{\bsl} \colon [0,1)^t \to \mathbb{C}$ as
\begin{align*}
&\trigarg{\bsl}{\bsy} = \exp(2 \pi \ii \bsl \cdot \bsy),
\end{align*}
where $\cdot$ denotes the usual Euclidean inner product.

Let now $s,t \in \NN$, $\alpha,\beta>1$ and let $\bsgamma^{(1)}, \bsgamma^{(2)}$ be two non-increasing sequences $\bsgamma^{(i)} = (\gamma_j^{(i)})_{j \geq 1}$ for $i \in \{1 , 2 \}$, where $0 < 
\gamma_j^{(i)} \leq 1$. We define two functions $\rho_{\alpha,\bsgamma^{(1)}}$ and
$r_{\beta,\bsgamma^{(2)}}$ as follows: For 
$\bsk = (k_1, \dots , k_s) \in \mathbb{N}_0^{s}$ and $\bsl = (l_1, \dots , l_t) \in \ZZ^t$ let
\begin{align*}
\rho_{\alpha,\bsgamma^{(1)}}(\bsk) = \prod_{j=1}^{s}\rho_{\alpha,\gamma_j^{(1)}}(k_j) \ \mbox{ and } \ r_{\beta,\bsgamma^{(2)}}(\bsl) = \prod_{j=1}^t r_{\beta,\gamma_j^{(2)}}(l_j),
\end{align*}
where
$$
\rho_{\alpha,\gamma_j^{(1)}}(k_j) = \begin{cases}
                    1  &   \text{ if } k_j = 0,\\
                    \gamma_j^{(1)} b^{-\alpha \left\lfloor \log_b(k_j) \right\rfloor} & \text{ if } k_j \neq 0,
                    \end{cases}
$$
and
$$
r_{\beta,\gamma_j^{(2)}}(l_j)= \begin{cases}
                    1  &   \text{ if } l_j = 0,\\
                    \gamma_j^{(2)} |l_j|^{-\beta} & \text{ if } l_j \neq 0.
                    \end{cases}
$$
With the help of these functions one can define so-called Walsh spaces \cite{DKPS, DP05} and Korobov spaces~\cite{DKS, LP, NW10}.

Here we define a hybrid function space as the tensor product of the Walsh and Korobov
spaces. The hybrid space $\mathcal{H}(K_{s,t,\alpha,\beta,\bsgamma})$, where $\bsgamma=(\bsgamma^{(1)},\bsgamma^{(2)})$, is the reproducing kernel Hilbert space with kernel function given by $K_{s,t,
\alpha,\beta,\bsgamma} : [0,1)^{s+t} \times [0,1)^{s+t} \rightarrow \mathbb{C}$,
\begin{align*}
K_{s,t,\alpha,\beta,\bsgamma}((\bsx, \bsy),(\bsx', \bsy')) = 
\sum_{\bsk \in \mathbb{N}_0^s}\sum_{\bsl \in \mathbb{Z}^t} \rho_{\alpha,\bsgamma^{(1)}}(\bsk)  r_{\beta,\bsgamma^{(2)}}(\bsl) \wal_{\bsk}(\bsx)\overline{\wal_{\bsk}(\bsx')}
\trigarg{\bsl}{\bsy}\overline{\trigarg{\bsl}{\bsy'}}
\end{align*}
and inner product
\begin{align*}
\left\langle f , g \right\rangle_{s,t,\alpha,\beta,\bsgamma} =
\sum_{\bsk \in \mathbb{N}_0^{s}}{\sum_{\bsl \in
\mathbb{Z}^{t}}\frac{1}{\rho_{\alpha,\bsgamma^{(1)}}(\bsk)} \frac{1}{r_{\beta,\bsgamma^{(2)}}(\bsl)}
\widehat{f}(\bsk,\bsl)\overline{\widehat{g}(\bsk,\bsl)}},
\end{align*}
with
\begin{align*}
\widehat{f}(\bsk,\bsl) = \int_{[0,1]^s} \int_{[0,1]^t}
{f(\bsx,
\bsy)\overline{\wal_{\bsk}(\bsx)\trigarg{\bsl}{\bsy}}}\,\mathrm{d}\bsx\mathrm{d}\bsy.
\end{align*}

The space $\mathcal{H}(K_{s,t,\alpha,\beta,\bsgamma})$ is the tensor product of a Walsh space and a Korobov space. 
If $s=0$, then we obtain the Korobov space, if $t=0$, then we obtain the Walsh space.

\begin{rem}\rm
For convenience we will in the following use the notation
$\int_{[0,1]^d} f(\bsx, \bsy) \, \mathrm{d} \bsx \mathrm{d} \bsy$, where $d=s+t$,
by which we mean
$\int_{[0,1]^{s}} \int_{[0,1]^{t}} f(\bsx,\bsy)\,\mathrm{d}\bsx\mathrm{d}\bsy$.
\end{rem}

The hybrid space $\mathcal{H}(K_{s,t,\alpha,\beta,\bsgamma})$ is the space of all absolutely convergent series $f$ of the form 
$$
f(\bsx,\bsy) = \sum_{(\bsk, \bsl) \in \mathbb{N}_0^{s} \times
\mathbb{Z}^{t}}{\widehat{f}(\bsk,
\bsl)\wal_{\bsk}(\bsx)\trigarg{\bsl}{\bsy}}\ \ \mbox{ with }\ \ \|f\|_{d,\alpha,\beta, \bsgamma} < \infty.$$
For further information on the space $\h$ we refer to \cite[Section~2.2]{KP15}.

\section{$\LL_2$-approximation}\label{secproblem}

Our goal is now to approximate the embedding from the hybrid space $\h$ to the space $\LL_2([0,1]^{s+t})$, i.e.,
$$\mathrm{EMB}_{s,t}: \h \rightarrow \LL_2([0,1]^{s+t}),\ \ \mathrm{EMB}_{s,t}(f)=f.$$ 
As already mentioned, it is enough to consider linear algorithms $A_{N,s,t}$ of the form
\begin{align}\label{eq:linalg}
A_{N,s,t}{(f)} = \sum_{k=1}^{N}a_k L_k(f),
\end{align}
with $a_k \in \LL_2([0,1]^{s+t})$ and continuous linear functionals $L_k$ on $\h$ from a permissible class of information $\Lambda$.  We consider two classes:

\begin{itemize}
\item $\Lambda=\Lambda^{\mathrm{all}}$ , the class of all continuous
linear functionals defined on $\h$.
Since  $\h$ is a Hilbert space, for every $L_k\in
\Lambda^{\mathrm{all}}$ there exists a 
function $f_k$ from $\h$ such that
$L_k(f)=\langle f,f_k\rangle_{d,\alpha,\beta,\bsgamma}$ for all $f\in
\h$.
\item $\Lambda=\Lambda^{\mathrm{std}}$, the class of
standard information consisting only of function evaluations. 
That is, $L_k\in\Lambda^{\mathrm{std}}$ if there exists
$(\bsx_k,\bsy_k)\in[0,1]^{s+t}$ such that $L_k(f)=f(\bsx_k,\bsy_k)$ for all $f\in \h$.
\end{itemize}

Since $\h$ is a reproducing kernel Hilbert space,
function evaluations are continuous linear functionals, and therefore
$\Lambda^{\mathrm{std}}\subseteq \Lambda^{\mathrm{all}}$. More
precisely,
$$
L_k(f)=f(\bsx_k,\bsy_k)=\langle f,K_{s,t,\alpha,\beta,\bsgamma}(\cdot,(\bsx_k,\bsy_k))\rangle_{s,t,\alpha,\beta,\bsgamma} $$
and
$$\|L_k\|=\|K_{s,t,\alpha,\beta,\bsgamma}\|_{s,t,\alpha,\beta,\bsgamma}=
K^{1/2}_{s,t,\alpha,\beta,\bsgamma}((\bsx_k,\bsy_k),(\bsx_k,\bsy_k)).
$$

The worst-case error in $\h$ of a linear algorithm as in \eqref{eq:linalg} is  
\begin{align*}
e^{\LL_2-\mathrm{app}}(A_{N,s,t}) = 
\sup_{\substack{f \in \h \\ \|f\|_{\h}\leq 1}}{\|\EMB_{s,t}{(f)} - A_{N,s,t}(f)\|_{\LL_2([0,1]^{s+t})}}.
\end{align*}
The $N^{{\rm th}}$ minimal worst-case error is given by
$$
e^{\LL_2-\mathrm{app},\Lambda}(N,s,t)=\inf_{A_{N,s,t}}e^{\mathrm{app}}(A_{N,s,t}),
$$
where the infimum is extended over all linear algorithms $A_{N,s,t}$ using information from the class $\Lambda$. The information complexity is given as 
\begin{align*}
N^{\LL_2-\mathrm{app},\Lambda}(\varepsilon,s,t) := \min\{ N : e^{\LL_2-\mathrm{app},\Lambda}(N,s,t) \leq \varepsilon\}.
\end{align*}
Since $\Lambda^{\rm std} \subseteq \Lambda^{\mathrm{all}}$, it follows that $N^{\LL_2-\mathrm{app},\Lambda^{\mathrm{all}}}(\varepsilon,s,t) \le N^{\LL_2-\mathrm{std},\Lambda^{\mathrm{std}}}
(\varepsilon,s,t)$.

We say that the $\LL_2$-approximation problem $\mathrm{EMB}=(\mathrm{EMB}_{s,t})_{s,t \ge 1}$ is:
\begin{itemize}
\item weakly tractable, if
\begin{align*}
\lim_{s +t + \varepsilon^{-1} \to \infty}{\frac{\log N^{\LL_2-\mathrm{app},\Lambda}(\varepsilon,s,t)}{s +t + \varepsilon^{-1}}} = 0;
\end{align*}
\item polynomially tractable, if we can find constants $C, \tau_1, \tau_2 \geq 0$ such that
\begin{align*}
N^{\LL_2-\mathrm{app},\Lambda}(\varepsilon,s,t) \leq C \varepsilon^{-\tau_1} (s+t)^{\tau_2} \text{ for all } 
\varepsilon \in (0,1) \text{ and all } s,t \in \mathbb{N};
\end{align*}
\item strongly polynomially tractable, if we can find constants $C, \tau_1\geq 0$ such that
\begin{align}\label{eqspt}
N^{\LL_2-\mathrm{app},\Lambda}(\varepsilon,s,t) \leq C \varepsilon^{-\tau_1} \text{ for all } \varepsilon \in (0,1) \text{ and all } s,t \in\mathbb{N}.
\end{align}
The infimum $\tau^{\ast}(\Lambda)$ of the real numbers $\tau_1$ such that \eqref{eqspt} 
holds is called the $\varepsilon$-exponent of strong polynomial tractability. 
\end{itemize}

For $\bsgamma=(\bsgamma^{(1)},\bsgamma^{(2)})$ we define the sum exponent  
\begin{equation}\label{def_sgamma} 
s_{\bsgamma}=\inf\left\{\kappa>0 \ : \ \sum_{j=1}^{\infty}(\gamma_j^{(1)})^{\kappa} < \infty \mbox{ and } \sum_{j=1}^{\infty}(\gamma_j^{(2)})^{\kappa} < \infty\right\}
\end{equation}
with the convention that $\inf \emptyset=\infty$.

Our main goal in this paper is to show the following theorem.

\begin{thm}\label{thm:main}
Consider the approximation problem $\mathrm{EMB}$. Then we have:
\begin{enumerate}
\item Strong polynomial tractability and polynomial tractability in the class $\Lambda^{\mathrm{all}}$ are equivalent, and they hold if and only if $s_{\bsgamma}< \infty$, where $s_{\bsgamma}$ is 
defined in \eqref{def_sgamma}. In this case the exponent of strong polynomial tractability is $\tau^{\ast}(\Lambda^{\mathrm{all}})=2 \max(s_{\bsgamma},\tfrac{1}{\alpha},\tfrac{1}{\beta})$. 

\item The problem is weakly tractable in the class $\Lambda^{\mathrm{all}}$ if and only if
\begin{align}\label{eq:necessaryWT}
\lim_{s+ t \to \infty}\frac{\sum_{j=1}^{s} \gamma_j^{(1)}+ \sum_{j=1}^{t}\gamma_j^{(2)}}{s+t} = 0.
\end{align}

\item The problem is strongly polynomially tractable in the class $\Lambda^{\mathrm{std}}$ if 
$$\sum_{j=1}^{\infty}\gamma_j^{(1)}< \infty \ \ \mbox{ and }\ \ \sum_{j=1}^{\infty}\gamma_j^{(2)}< \infty.$$ 
The exponent of strong polynomial tractability in the class $\Lambda^{\mathrm{std}}$ 
satisfies 
$$\tau^{\ast}(\Lambda^{\mathrm{std}})\in [2 \max(\tfrac{1}{\alpha},\tfrac{1}{\beta},s_{\bsgamma}),4+2 \max(\tfrac{1}{\alpha},\tfrac{1}{\beta},s_{\bsgamma})].$$ 

\item The problem is polynomially tractable in the class $\Lambda^{\mathrm{std}}$ if $$\limsup_{s\To\infty}\frac{\sum_{j=1}^{s}\gamma_j^{(1)}}{\log s}< \infty \ \ \mbox{ and }\ \ 
\limsup_{t\To\infty}\frac{\sum_{j=1}^{t}\gamma_j^{(2)}}{\log t}< \infty.$$

\item The problem is weakly tractable in the class $\Lambda^{\mathrm{std}}$ if and only if
$$\lim_{s+ t \to \infty}\frac{\sum_{j=1}^{s} \gamma_j^{(1)}+ \sum_{j=1}^{t}\gamma_j^{(2)}}{s+t} = 0.$$
\end{enumerate}
\end{thm}

\begin{rem}\rm
 Since it can easily be verified that integration in $\h$ is not harder than approximation, the last item in Theorem \ref{thm:main} implies that the 
 sufficient condition for weak tractability of integration shown in \cite{KP15} is also necessary. 
\end{rem}

\section{Proof of Theorem~\ref{thm:main}}\label{sec:proof}

We recall that strong polynomial tractability implies polynomial tractability which in turn implies weak tractability. 
Furthermore, all sufficient conditions for the class 
$\Lambda^{\mathrm{std}}$ are also sufficient 
for the class $\Lambda^{\mathrm{all}}$ with $\tau^*(\Lambda^{\mathrm{all}}) \le \tau^*(\Lambda^{\mathrm{std}})$ in the case of strong polynomial tractability. 
All necessary  conditions for the class $\Lambda^{\mathrm{all}}$ are also necessary for the class $\Lambda^{\mathrm{std}}$.

\subsection{Proof of Item 1} 

In order to give a necessary and sufficient condition for strong polynomial tractability for $\Lambda^{\rm all}$  
we use a criterion from \cite[Section~5.1]{NW08}. Let us consider the self-adjoint operator 
$W_{s,t}:=\mathrm{EMB}_{s,t}^* \mathrm{EMB}_{s,t}: \h \rightarrow \h$, which in our case is 
given by 
$$ 
W_{s,t} f=\sum_{(\bsk,\bsl) \in \NN_0^{s}\times\ZZ^{t}}\rho_{\alpha,\bsgamma^{(1)}}(\bsk) r_{\beta,\bsgamma^{(2)}}(\bsl)\widehat{f}(\bsk,\bsl)\wal_{\bsk}(\bsx)\trigarg{\bsl}{\bsy}.
$$
The eigenvalues are then given by the collection of the numbers
$$\rho_{\alpha,\bsgamma^{(1)}}(\bsk) r_{\beta,\bsgamma^{(2)}}(\bsl) \ \ \ \mbox{ for $(\bsk,\bsl) \in \NN_0^{s}\times\ZZ^{t}$}.$$ 
Furthermore, the largest eigenvalue is
$\rho_{\alpha,\bsgamma^{(1)}}(\bszero) r_{\beta,\bsgamma^{(2)}}(\bszero)=1$. 

From \cite[Theorem 5.2]{NW08} we find that the problem $\mathrm{EMB}$ is polynomially tractable for $\Lambda^{\rm all}$ if and only if there exist $\nu>0$ and $q\ge 0$ such that 
\begin{equation}\label{critNW08} 
\sup_{s,t}\left(\sum_{(\bsk,\bsl) \in \NN_0^{s}\times\ZZ^{t}} (\rho_{\alpha,\bsgamma^{(1)}}(\bsk) r_{\beta,\bsgamma^{(2)}}(\bsl))^\nu \right)^{1/\nu} (s+t)^{-q}< \infty.
\end{equation}
Furthermore, we have strong polynomial tractability if and only if \eqref{critNW08} holds with $q=0$.

It is easy to check that we require $\nu>\max(\tfrac{1}{\alpha},\tfrac{1}{\beta})$ in order for \eqref{critNW08} to hold with $q=0$. Let us now assume that $\nu$ is indeed bigger than $\max(\tfrac{1}
{\alpha},\tfrac{1}{\beta})$. For the sum in \eqref{critNW08} we have
\begin{eqnarray}\label{eqnecessary}
 \sum_{(\bsk,\bsl) \in \NN_0^{s}\times\ZZ^{t}} (\rho_{\alpha,\bsgamma^{(1)}}(\bsk) r_{\beta,\bsgamma^{(2)}}(\bsl))^\nu
 &=&\prod_{j=1}^{s}\left(1+ (\gamma_j^{(1)})^\nu \mu (\alpha\nu)\right)
    \prod_{j=1}^{t}\left(1+ (\gamma_j^{(2)})^\nu 2\zeta (\beta\nu)\right),
\end{eqnarray}
where $\mu(x)=\frac{b^{x}(b-1)}{b^{x}-b}$ for $x>1$ and $\zeta(x)$ is the Riemann zeta function.

Now, using arguments outlined in \cite{SW01} (see also \cite[Section~4.5]{LP}), it is easy to see that the existence of some $\nu>\max(\tfrac{1}{\alpha},\tfrac{1}{\beta})$ with  
\begin{equation*}
\sum_{j=1}^{\infty}(\gamma_j^{(1)})^\nu< \infty \quad\mbox{and}\quad \sum_{j=1}^{\infty}(\gamma_j^{(2)})^\nu< \infty
\end{equation*}
 is a necessary and sufficient condition for \eqref{critNW08} with $q=0$ and therefore for strong polynomial tractability of the problem $\mathrm{EMB}$.

Again according to \cite[Theorem 5.2]{NW08}, the exponent of strong polynomial tractability is $2 \max(\tfrac{1}{\alpha},\tfrac{1}{\beta},s_{\bsgamma})$, where $s_{\bsgamma}$ is defined in 
\eqref{def_sgamma}.

It remains to show the equivalence of strong polynomial and polynomial tractability. Of course, it suffices to show that polynomial tractability implies strong polynomial tractability. So assume that 
the problem $\mathrm{EMB}$ is polynomially tractable for the class $\Lambda^{\mathrm{all}}$. Then we obtain polynomial tractability also for the embedding problem in the pure Walsh space $\mathcal{H}
(K_{s,0,\alpha,\beta,\bsgamma})$ and in the pure Korobov space $\mathcal{H}(K_{0,t,\alpha,\beta, \bsgamma})$. According to \cite[Theorem~2]{WW99} this is equivalent to strong polynomial tractability  
for the embedding problem in the pure Walsh space $\mathcal{H}(K_{s,0,\alpha,\beta,\bsgamma})$ and in the pure Korobov space $\mathcal{H}(K_{0,t,\alpha,\beta, \bsgamma})$. According to \cite{DKK08} and 
\cite{KSW06} this implies the existence of $\nu_1>0$ such that $\sum_{j \ge 1} (\gamma_j^{(1)})^{\nu_1} < \infty$ and the existence of $\nu_2>0$ such that $\sum_{j \ge 1} (\gamma_j^{(2)})^{\nu_2} < 
\infty$. Hence we have $s_{\bsgamma}< \infty$ and this in turn implies strong polynomial tractability for the class $\Lambda^{\mathrm{all}}$, as shown above. This completes the proof of Item 1.

\subsection{Proof of Item 2}

Sufficiency of Condition \eqref{eq:necessaryWT} follows by Item 5 of the Theorem which we show in the next section. 

For showing necessity of Condition \eqref{eq:necessaryWT}, we use \cite[Theorem~5.3]{NW08} in the following. To keep notation simple, we shall frequently write again
$d$ instead of $s+t$. Theorem~5.3 in \cite{NW08} states that our approximation problem is weakly tractable for $\Lambda^{\rm all}$ if and only if 
\begin{itemize}
	\item $\lim\limits_{j \to \infty}\lambda_{d,j}\log^2j = 0$ for all $d \in \NN$ and
	\item there exists some function $f \colon (0,\frac{1}{2}] \to \NN$ such that
	\begin{align}\label{eq:supCondition}
	\sup\limits_{\eta \in (0, \frac{1}{2}]} \frac{1}{\eta^2} \sup\limits_{d \geq f(\eta)} \sup\limits_{j \geq \left\lceil \exp(d\sqrt{\eta}) \right\rceil + 1} 
	\lambda_{d,j}\log^2j < \infty,
	\end{align}
\end{itemize}
where $\lambda_{d,j}=\lambda_{s+t,j}$ denotes the $j^{{\rm th}}$ eigenvalue of $W_{s,t}$ in non-increasing order.

Let us now assume that the approximation problem is weakly tractable for $\Lambda^{\rm all}$. This then in particular implies that 
\begin{equation}\label{eqcondeigen}
\lim\limits_{j \to \infty}\lambda_{d,j}\log^2j = 0\quad\mbox{for all}\quad d \in \NN.
\end{equation}
We are now going to show that \eqref{eqcondeigen} implies \eqref{eq:necessaryWT}. To this end, recall that the eigenvalues of $W_{s,t}$ are of the form 
$$\rho_{\alpha,\bsgamma^{(1)}}(\bsk) r_{\beta,\bsgamma^{(2)}}(\bsl) \ \ \ \mbox{ for $(\bsk,\bsl) \in \NN_0^{s}\times\ZZ^{t}$}.$$
Note that we have $\lambda_{d,1}=1$; furthermore, note that $\rho_{\alpha,\gamma_j^{(1)}}(1)=\gamma_j^{(1)}$ for any $j\in\NN$, and 
$r_{\beta,\gamma_i^{(2)}}(1)=\gamma_i^{(2)}$ for any $i\in\NN$. Hence, by choosing all components of $(\bsk,\bsl)\in \NN_0^{s}\times\ZZ^{t}$ but one 
equal to zero, and the remaining equal to one, we see that 
$$\gamma_1^{(1)},\ldots,\gamma_s^{(1)}\quad\mbox{and}\quad \gamma_1^{(2)},\ldots,\gamma_t^{(2)}$$ 
are eigenvalues of $W_{s,t}$. Consequently,
$$\sum_{j=1}^s \gamma_j^{(1)} + \sum_{j=1}^t \gamma_j^{(2)}\le \sum_{j=1}^d \lambda_{d,j},$$
and hence
$$\lim_{s+ t \to \infty}\frac{\sum_{j=1}^{s} \gamma_j^{(1)}+ \sum_{j=1}^{t}\gamma_j^{(2)}}{s+t} \le 
  \lim_{d \to \infty}\frac{\sum_{j=1}^{d} \lambda_{d,j}}{d}.$$
However, due to \eqref{eqcondeigen}, it follows that the latter limit is 0, which shows that indeed \eqref{eq:necessaryWT} holds.

\subsection{Proof of Items 3--5}

Any $f \in \h$ can be displayed as
$$
f(\bsx ,\bsy) = \sum_{(\bsk, \bsl) \in \mathbb{N}_0^{s}
\times \mathbb{Z}^{t}}{\widehat{f}(\bsk,
\bsl)\walarg{\bsk}{\bsx}\trigarg{\bsl}{\bsy}}.
$$

In order to approximate $\widehat{f}(\bsk,\bsl)$, we are going to use quasi-Monte Carlo algorithms 
based on classical and on polynomial lattice point sets.
\paragraph{Classical lattice point sets.}
For $N \in \NN$ and  $\boldz=(z_1,\ldots,z_t) \in
Z_N^t$, where $Z_N:=\{z\in\{1,\ldots,N-1\}: \gcd (z,N)=1\}$, the lattice point set $\{\boldq_v\}_{v=0}^{N-1}$
with generating vector $\boldz$ is defined by
\[\boldq_v=\left(\left\{\frac{v z_1}{N}\right\},\ldots,\left\{\frac{v z_{t}}{N}\right\}\right)
 \; \mbox{ for all }\; 0\le v\le N-1.\] Here $\{\cdot\}$ denotes the fractional part
of a real number.

\paragraph{Polynomial lattice point sets.}
Let $\FF_b$ be the finite field of prime order $b$, $\FF_b[x]$ be the set of polynomials over
$\FF_b$, and let $\FF_b ((x^{-1}))$ be the field of formal Laurent
series over $\FF_b$. The latter contains the field of rational functions
as a subfield. Given $m\in \NN$, set $G_{b,m}:=\{a \in \FF_b[x]\, : \, \deg(a)<m\}$ and define a mapping
$\nu_m:\FF_b ((x^{-1}))\To [0,1)$ by
\[\nu_m\left(\sum_{l=z}^{\infty}t_l x^{-l}\right):=\sum_{l=\max (1,z)}^{m}t_l b^{-l}.\]
Let $f\in\FF_b [x]$ with $\deg(f)=m$ and $\boldg=(g_1,\ldots,g_{s})\in\FF_b [x]^{s}$. 
The polynomial lattice point set $(\boldp_v)_{v \in G_{b,m}}$ 
with generating vector $\boldg$  
is defined by   
\[\boldp_v:=\left(\nu_m\left(\frac{v(x)g_1 (x)}{f(x)}\right),
\ldots,\nu_m\left(\frac{v(x)g_{s} (x)}{f(x)}\right)\right) \; \mbox{ for all }\; v \in G_{b,m}.\]
Note that we can associate the polynomial $v(x)=\sum_{r=0}^{m-1} v_r x^r \in G_{b,m}$ with the integer 
$v=\sum_{r=0}^{m-1} v_r b^r$, where, with a slight abuse of notation, the element $v_r\in\FF_b$ is 
associated with the integer $v_r\in\{0,1\ldots,b-1\}$. In this way we can index the points of a 
polynomial lattice point set by integers ranging from $0$ to $b^m-1$. 

\medskip

Now suppose that $N$ is of the form $b^m$ for some $m\in\NN$, 
and let $\mathrm{PL} = \{ \boldp_0, \dots , \boldp_{N-1}\} \subseteq
[0,1)^s$ be a polynomial lattice point set and $\mathrm{L} = \{
\boldq_0, \dots , \boldq_{N-1}\} \subseteq [0,1)^t$ be a
lattice point set. We consider the point set $\mathrm{(PL,L)} = \{ (\boldp, \boldq)_v=(\boldp_v,\boldq_v) \colon v = 0 , \dots, N-1\}$.

For $M \ge 1$ define the set
\begin{equation}\label{eq:aall}
\mathcal{A}_M = \{ (\bsk, \bsl) \in \mathbb{N}_0^{s} \times
\mathbb{Z}^{t} :  ( \rho_{\alpha,\bsgamma^{(1)}}(\bsk) )^{-1}
(r_{\beta,\bsgamma^{(2)}}(\bsl))^{-1} \leq M \}.
\end{equation}

In order to approximate the embedding $\mathrm{EMB}_{s,t}(f)=f$ for $f \in \h$ we use the algorithm  
\begin{align}\label{eq:algUsed}
\linalg{\bsx, \bsy} = \sum_{(\bsk,\bsl) \in \mathcal{A}_M}{\left(
\frac{1}{N}\sum_{v=0}^{N-1}{f((\boldp,
\boldq)_v)\overline{\walarg{\bsk}{\boldp_v}\trigarg{\bsl}{\boldq_v}}}\right)\walarg{\bsk}{\bsx}\trigarg{\bsl}{\bsy}}.
\end{align}
It can easily be checked that $A_{N,s,t,M}$ is a linear algorithm of the form \eqref{eq:linalg} with
\begin{align*}
a_v(\bsx,\bsy) = \frac{1}{N} \sum_{(\bsk,\bsl) \in
\mathcal{A}_M}{\wal_{\bsk}(\bsx \ominus
\boldp_v)\trigarg{\bsl}{\bsy - \boldq_v}}\ \ \text{ and }\ \ L_v(f)= f((\boldp, \boldq)_v),\ 0\le v\le N-1.
\end{align*}
The error of approximation for given $f \in \h$ is then
\begin{align}\label{eq:errorapprox}
(f - A_{N,s,t,M}(f))(\bsx, \bsy) &= 
\sum_{(\bsk, \bsl) \notin \mathcal{A}_M}{\widehat{f}(\bsk, \bsl)\wal_{\bsk}(\bsx)\trigarg{\bsl}{\bsy}} \nonumber \\
&\quad+ \sum_{(\bsk, \bsl) \in \mathcal{A}_M}{\left( \widehat{f}(\bsk, \bsl) - 
\frac{1}{N}\sum_{v=0}^{N-1}{f((\boldp, \boldq)_v)\overline{\wal_{\bsk}(\boldp_v)\trigarg{\bsl}{\boldq_v}}}\right)
\walarg{\bsk}{\bsx}\trigarg{\bsl}{\bsy}}.
\end{align}
We use \eqref{eq:errorapprox} and Parseval's identity to obtain
\begin{eqnarray*}
\| \EMB_{s,t}(f) - A_{N,s,t,M}(f) \|_{\LL_2([0,1]^{s + t})}^2 = S_1 + S_2,
\end{eqnarray*}
where 
\[
 S_1:=\sum_{(\bsk, \bsl) \notin \mathcal{A}_M} | \widehat{f}(\bsk, \bsl)|^2, 
\]
and
\[
 S_2:=\sum_{(\bsk, \bsl) \in \mathcal{A}_M} \left| \int_{[0,1]^{s + t}} f_{(\bsk, \bsl)}(\bsx, \bsy) \, \mathrm{d}\bsx\mathrm{d} \bsy - 
\frac{1}{N} \sum_{v=0}^{N-1} f_{(\bsk, \bsl)}((\boldp, \boldq)_v)  \right|^2,
\]
with
\begin{align*}
f_{(\bsk, \bsl)}(\bsx, \bsy) := f(\bsx, \bsy)
\overline{\walarg{\bsk}{\bsx} \trigarg{\bsl}{\bsy} }.
\end{align*}

From the definition of $\mathcal{A}_M$ it follows easily that
\begin{align*}
S_1 < \frac{1}{M} \| f \|_{\h}^2.
\end{align*}

Let us now consider $S_2$. The term in-between the absolute value signs in $S_2$ is 
the integration error of the QMC rule using the nodes $\mathrm{(PL,L)}$ for the function $f_{(\bsk, \bsl)}(\bsx, \bsy)$. Since the product of two Walsh functions 
is again a Walsh function, and the analogue is true for trigonometric functions, 
it can easily be verified that $f_{(\bsk, \bsl)}\in \h$.
Hence we can bound $S_2$ by
\begin{align*}
S_2 &\leq 
(e^{\mathrm{int}}(\mathrm{PL,L}))^2 \sum_{(\bsk, \bsl) \in \mathcal{A}_M}  \| f_{(\bsk, \bsl)} \|_{\h}^2,
\end{align*}
where $e^{\mathrm{int}}(\mathrm{PL,L})$ is the worst-case integration error in $\h$ of the QMC rule based on the nodes $\mathrm{(PL,L)}$, i.e.,
\begin{align*}
e^{\mathrm{int}}(\mathrm{PL,L})= \sup_{\substack{f \in \h \\ \|f\|_{\h}\le 1}} \left| \int_{[0,1]^{s + t}} f(\bsx, \bsy) \, \mathrm{d}\bsx\mathrm{d} \bsy - 
\frac{1}{N} \sum_{v=0}^{N-1} f((\boldp, \boldq)_v)  \right|.
\end{align*}

From \cite[Theorem~3]{KP15} it then follows that
\begin{align}\label{eq:S_2}
S_2 
&\leq  \frac{2}{N} 
\left( \prod_{j=1}^{s}(1 + \gamma_j^{(1)} 2\mu(\alpha)) \right) 
\left( \prod_{j=1}^{t}(1 + \gamma_j^{(2)} 4\zeta(\beta)) \right) \sum_{(\bsk, \bsl) \in \mathcal{A}_M} \| f_{(\bsk, \bsl)} \|_{\h}^2.
\end{align}

Next we find an estimate for $\| f_{(\bsk, \bsl)}\|_{\h}^2$ for $(\bsk, \bsl) \in \mathcal{A}_M$. 
From the easily seen fact that $\widehat{f}_{(\bsk, \bsl)}{(\bsh, \bsm)} = \widehat{f}(\bsk \oplus 
\bsh, \bsl + \bsm)$ we obtain
\begin{align*}
\| f_{(\bsk, \bsl)} \|_{\h}^2 &= \sum_{\bsh \in \mathbb{N}_0^{s}} \sum_{\bsm \in \mathbb{Z}^{t}} \frac{| \widehat{f}(\bsk \oplus \bsh, \bsl + \bsm) |^2}{\rho_{\alpha,\bsgamma^{(1)}}(\bsh) r_{\beta,
\bsgamma^{(2)}}(\bsm)} \\
&= \sum_{\bsh \in \mathbb{N}_0^{s}} \sum_{\bsm \in \mathbb{Z}^{t}}  \frac{| \widehat{f}(\bsk \oplus \bsh, \bsl + \bsm) |^2}{\rho_{\alpha,\bsgamma^{(1)}}(\bsk \oplus \bsh)r_{\beta,\bsgamma^{(2)}}(\bsl + 
\bsm)} \frac{ \rho_{\alpha,\bsgamma^{(1)}}(\bsk \oplus \bsh)  r_{\beta,\bsgamma^{(2)}}(\bsl + \bsm)}{\rho_{\alpha,\bsgamma^{(1)}}(\bsh)  r_{\beta,\bsgamma^{(2)}}(\bsm)}.
\end{align*}

Combining results from \cite{DKK08} and \cite{KSW06} we find
\begin{align*}
\frac{ \rho_{\alpha,\bsgamma^{(1)}}(\bsk \oplus \bsh)  r_{\beta,\bsgamma^{(2)}}(\bsl + \bsm)}{\rho_{\alpha,\bsgamma^{(1)}}(\bsh)  r_{\beta,\bsgamma^{(2)}}(\bsm)} \le \frac{1}{\rho_{\alpha,
\bsgamma^{(1)}}(\bsk)  r_{\beta,\bsgamma^{(2)}}(\bsl)} 
\prod_{j=1}^{t} \max(1, 2^{\beta}\gamma_j^{(2)}) \le M \prod_{j=1}^{t} \max(1,2^{\beta}\gamma_j^{(2)}),
\end{align*}
and hence, after applying an index shift,
\begin{align*}
\|f_{(\bsk, \bsl)}\|_{\h}^2 &\leq M
\prod_{j=1}^{t} \max(1, 2^{\beta}\gamma_j^{(2)}) \sum_{\bsh \in \mathbb{N}_0^{s}} \sum_{\bsm \in \mathbb{Z}^{t}}
\frac{ |\widehat{f}(\bsk \oplus \bsh, \bsl + \bsm)|^2 }{\rho_{\alpha,\bsgamma^{(1)}}(\bsk \oplus \bsh) r_{\beta,\bsgamma^{(2)}}(\bsl + \bsm)}\\
&= M \| f \|_{\h}^2
\prod_{j=1}^{t} \max(1, 2^{\beta}\gamma_j^{(2)}).
\end{align*}

Plugging this into \eqref{eq:S_2} we obtain
\begin{align}\label{eq:S_2zwei}
S_2 &\leq \frac{2}{N} \left( \prod_{j=1}^{s}(1 + \gamma_j^{(1)}2\mu(\alpha)) \right) 
\left( \prod_{j=1}^{t}(1 + \gamma_j^{(2)}4\zeta(\beta)) \right) \| f\|_{\h}^2 
 M |\mathcal{A}_M| \prod_{j=1}^{t} \max(1, 2^{\beta}\gamma_j^{(2)}).
\end{align}

Next we study the cardinality of the set $\mathcal{A}_M$.

\begin{lem}\label{lem:aall}
Let $\theta=\min(\alpha,\beta)$. For arbitrary $\kappa > 1/\theta=\max(\tfrac{1}{\alpha},\tfrac{1}{\beta})$ we have
$$|\mathcal{A}_M| \le M^{\kappa} \prod_{j=1}^s \left(1+2 \zeta(\theta \kappa) (b^{\alpha} \gamma_j^{(1)})^{\kappa}\right) 
\prod_{j=1}^t \left(1+2 \zeta(\theta \kappa) (\gamma_j^{(2)})^{\kappa}\right).$$ 
\end{lem}

\begin{proof}
For $k \in \NN$ we have 
$$\frac{1}{\rho_{\alpha,\gamma}(k)}=\frac{b^{\alpha \lfloor \log_b k\rfloor}}{\gamma} \ge \frac{b^{\alpha (-1+\log_b k)}}{\gamma} = \frac{k^{\alpha}}{\gamma b^{\alpha}}=\frac{1}
{r_{\alpha,\gamma b^{\alpha}}(k)}.$$  

Then we have
\begin{align*}
\mathcal{A}_M & =  \left\{ (\bsk, \bsl) \in \mathbb{N}_0^s
\times \mathbb{Z}^t : \frac{1}{\rho_{\alpha,\bsgamma^{(1)}}(\bsk)} \frac{1}{r_{\beta,\bsgamma^{(2)}}(\bsl)} \leq M \right\}\\
& \subseteq \left\{ (\bsk, \bsl) \in \NN_0^s
\times \mathbb{Z}^{t} \colon \frac{1}{r_{\alpha,\bsgamma^{(1)} b^{\alpha}}(\bsk)} \frac{1}{r_{\beta,\bsgamma^{(2)}}(\bsl)} 
\leq M \right\}\\
& \subseteq \left\{ (\bsk, \bsl) \in \ZZ^s
\times \mathbb{Z}^{t} \colon \frac{1}{r_{\theta,\bsgamma^{(1)} b^{\alpha}}(\bsk)} \frac{1}{r_{\theta,\bsgamma^{(2)}}(\bsl)} 
\leq M \right\}
\end{align*}
from which the result follows immediately from \cite[Lemma~1]{KSW06}.
\end{proof}

Considering Lemma~\ref{lem:aall}, for any $\kappa > 1/\nu$ we obtain
$$S_2 \le c_{s,t,\alpha,\beta,\bsgamma,\kappa} \frac{M^{1+\kappa}}{N} \| f \|_{\h}^2,$$
where
\begin{align}\label{eq:constOhne}
c_{s,t,\alpha,\beta,\bsgamma,\kappa} &:= 2 \left( \prod_{j=1}^{s}(1 + \gamma_j^{(1)}2\mu(\alpha)) \right) 
\left( \prod_{j=1}^{t}(1 + \gamma_j^{(2)}4\zeta(\beta)) \right) 
\prod_{j=1}^{t} \max(1, 2^{\beta}\gamma_j^{(2)})  \nonumber\\
&\quad\times \prod_{j=1}^{s} \left(1+2 \zeta(\theta \kappa) (b^{\alpha} \gamma_j^{(1)})^{\kappa}\right) 
\prod_{j=1}^{t} \left(1+2 \zeta(\theta \kappa) (\gamma_j^{(2)})^{\kappa}\right).
\end{align}

Summing up we have
\begin{align*}
\left\lVert \EMB_{s,t}(f) - A_{N,s,t,M}(f)
\right\rVert^2_{\LL_2{([0,1]^{s +t})}} \leq \left( \frac{1}{M} +
c_{s,t,\alpha,\beta,\bsgamma,\kappa} \frac{M^{1+\kappa}}{N} \right) \| f \|_{\h}^2.
\end{align*}
Choosing $M=M(N) = (N/c_{s,t,\alpha,\beta,\bsgamma,\kappa})^{1/(2+\kappa)}$ and taking the square root we obtain 
the following proposition and its corollary, which then concludes the proof of Theorem \ref{thm:main}.
\begin{prop}\label{thm_err_bd}
Let $\kappa>1/\min(\alpha,\beta)$ and let $c_{s,t,\alpha,\beta,\bsgamma,\kappa}$ be defined as in \eqref{eq:constOhne}. The worst-case error of the algorithm $A_{N,s,t,M}$ as defined in 
\eqref{eq:algUsed} 
using a point set $\mathrm{(PL,L)}$ constructed by \cite[Algorithm 1]{KP15} and with $M= (N/c_{s,t,\alpha,\beta,\bsgamma,\kappa})^{1/(2+\kappa)}$
satisfies 
\begin{align*}
e^{\LL_2-\mathrm{app}}(A_{N,s,t,M}) \leq \sqrt{2} \left( \frac{c_{s,t,\alpha,\beta,\bsgamma,\kappa}}{N}\right)^{\frac{1}{4+2\kappa}}.
\end{align*}
\end{prop}

\begin{cor}\label{thm:sufftract}
Consider the approximation problem $\mathrm{EMB}$ with information from the class $\Lambda^{{\rm std}}$.
\begin{itemize}
\item If $\sum_{j=1}^{\infty}\gamma_j^{(1)}< \infty$ and $\sum_{j=1}^{\infty}\gamma_j^{(2)}< \infty$, then $\mathrm{EMB}$ is strongly polynomially tractable with $\varepsilon$-exponent at most
$4+2\max(s_{\bsgamma}, \tfrac{1}{\alpha},\tfrac{1}{\beta})$; 
\item if $\limsup_{s \to \infty}\sum_{j=1}^{s}\frac{\gamma_j^{(1)}}{\log{(s+1)}}
< \infty$ and $\limsup_{t \to
\infty}\sum_{j=1}^{t}\frac{\gamma_j^{(2)}}{\log{(t+1)}}
< \infty$, then $\mathrm{EMB}$ is polynomially tractable; 
\item if $\lim_{s+ t \to \infty} \frac{\sum_{j=1}^{s}
\gamma_j^{(1)} +\sum_{j=1}^{t}\gamma_j^{(2)}}{s+t} = 0$, then $\mathrm{EMB}$ is weakly tractable.
\end{itemize}
\end{cor}

\begin{proof}
Employing Proposition~\ref{thm_err_bd}, the result follows by the same arguments as used in \cite[Section~5.2]{KP15}. 
We only show the first item:  Let $\kappa=1$. Since $\log(1+x)\le x$ we obtain 
\begin{eqnarray*}
c_{s,t,\alpha,\beta,\bsgamma,1} \le 2 \exp\left(u_1(\alpha,\beta) \sum_{j=1}^{s} \gamma_j^{(1)} +u_2(\alpha,\beta) \sum_{j=1}^{s} \gamma_j^{(2)} \right) \le c_{\infty,\infty,\alpha,\beta,\bsgamma,1} < 
\infty,
\end{eqnarray*}
where $u_1(\alpha,\beta)=2 \mu(\alpha)+2 \zeta(\theta) b^{\alpha}$ and $u_2(\alpha,\beta)=4 \zeta(\beta)+2^{\beta}+2\zeta(\theta)$.
Then Proposition~\ref{thm_err_bd} with $\kappa=1$ implies that 
$$e^{\LL_2-\mathrm{app}}(N,s+t) \le \sqrt{2} \left(\frac{c_{\infty,\infty,\alpha,\beta,\bsgamma,1}}{N}\right)^{1/6}.$$
Recall that $N$ is of the form $b^m$. Now, for $\varepsilon>0$ choose $m \in \NN$ such that 
$b^{m-1} < \lceil 8  c_{\infty,\infty,\alpha,\beta,\bsgamma,1} \varepsilon^{-6}\rceil \linebreak =: N' \le b^m$. Then we have 
$e^{\LL_2-\mathrm{app}}(b^m,s+t)\le \varepsilon$ 
and hence $$N^{\LL_2-\mathrm{app},\Lambda^{\mathrm{std}}}(\varepsilon,s+t) \le b^m < b N' =b \lceil 8  c_{\infty,\infty,\alpha,\beta,\bsgamma,1} \varepsilon^{-6}\rceil.$$ This is strong polynomial 
tractability. The result for the $\varepsilon$-exponent can be shown easily by similar arguments.
\end{proof}

\section*{Acknowledgments}
P. Kritzer would like to thank R.F. Tichy and S. Thonhauser for their hospitality 
during his stay at Graz University of Technology, where parts of this paper were written.

\begin{small}
\noindent\textbf{Authors' addresses:}\\

\medskip

\noindent Peter Kritzer and Helene Laimer\\
Johann Radon Institute for Computational and Applied Mathematics (RICAM)\\
Austrian Academy of Sciences\\
Altenbergerstr.~69, 4040 Linz, Austria\\
E-mail: \texttt{peter.kritzer@oeaw.ac.at}, \texttt{helene.laimer@oeaw.ac.at}

\medskip

\noindent Friedrich Pillichshammer\\
Institut f\"{u}r Finanzmathematik und Angewandte Zahlentheorie\\
Johannes Kepler Universit\"{a}t Linz\\
Altenbergerstr.~69, 4040 Linz, Austria\\
E-mail: \texttt{friedrich.pillichshammer@jku.at}

\end{small}
\end{document}